\documentclass{amsart}
\newtheorem{theorem}{Theorem}[section]

\newtheorem{corollary}[theorem]{Corollary}
\newtheorem{proposition}[theorem]{Proposition}

\theoremstyle{definition}

\theoremstyle{remark}
\newtheorem{remark}[theorem]{Remark}

\numberwithin{equation}{section}

 \def \l{{\lambda}}
\def \a{{\alpha}}
\def \b{{\beta}}
\def \s{{\sigma}}

\def \O{{\Omega}}
\def \d{{\delta}}
\def \B{{\mathcal B}}

\def \N{{\bf N}}
\def \m{{\mu}}
 
\def \D{{\Delta}}

 \def \P{{\bf P}}

\def \qq{{\qquad}}
 \def \p{{\varphi}}
\def \noi{{\noindent}}

\def \e{{\varepsilon}}
 \def \t{{\theta}}

 at 10,4 pt

 at 10,8 pt
 \def \lt {{\hbox{\vrule height 6pt depth 2pt width 0,7pt}\kern 1,1pt}}
\def \rt {{\kern 0,9pt\hbox{\vrule height 6pt depth 2pt width 0,7pt}}}

\font\gsec= cmbx10 at 9 pt \scrollmode
  
  \scrollmode
\font\gsec= cmb10 at 10 pt

\def\P{{\mathbb P}}

\def\N{{\mathbb N}}

at  8 pt
\scrollmode

\font\gsec= cmb10 at 10 pt
   at 10,4  pt at  8,2 pt
  \scrollmode
\font\gsec= cmb10 at 9    pt
 at 8,2  pt
 
\begin{document}

\title[{ Zerofree Regions of the Riemann Zeta-Function}
]{Local Suprema of Dirichlet Polynomials and Zerofree Regions of the Riemann Zeta-Function}

\author{Michel J.\,G. Weber}  
 \address{Michel Weber: IRMA, 10  rue du G\'en\'eral Zimmer, 67084
 Strasbourg Cedex, France}
  \email{michel.weber@math.unistra.fr \!;  \! m.j.g.weber@gmail.com}

 %%%%%%%%%%%%%%%%%%%%%%%%%%%%%%%%%%%%%%%%%%%%%%%%%%%%%%%%%%%%%%%%%%%%%%%%%%%%%%%%%%%%%%%%%%%%%%

\begin{abstract}     A new   zerofree region  of  the Riemann Zeta-function $\zeta$ is identified by using Tur\'an's
localization criterion     linking   zeros  of    $\zeta$ with  
uniform  local suprema of sets of
Dirichlet polynomials expanded over the primes. The proof is based on a randomization argument. An estimate  for  local extrema    for some
finite families   of 
  shifted  Dirichlet polynomials, 
% associated to linearly independent   sequences of reals
 is     established by  preliminary considering their local
increment properties, by means of Montgomery-Vaughan's variant of  Hilbert's inequality.  A   covering
argument  combined with Tur\'an's
localization criterion  allows to conclude.  
 \end{abstract}
 %%%%%%%%%%%%%%%%%%%%%%%%%%%%%%%%%%%%%%%%%%%%%%%%%%%%%%%%%%%%%%%%%%%%%%%%%%%%%%%%%%%%%%%%%%%%%%%

\maketitle

%%%%%%%%%%%%%%%%%%%%%%%%%%%%%%%%%%%%%%%%%%%%%%%%%%%%%%%%%%%%%%%%%%%%%%%%%%%%%%%%%%%%%%%%%%%%%%%%%
%%%%%%%%%%%%%%%%%%%%%%%%%%%%%%%%%%%%%%%%%%%%%%%%%%%%%%%%%%%%%%%%%%%%%%%%%%%%%%%%%%%%%%%%%%%%%%%%%
%%%%%%%%%%%%%%%%%%%%%%%%%%%%%%%%%%%%%%%%%%%%%%%%%%%%%%%%%%%%%%%%%%%%%%%%%%%%%%%%%%%%%%%%%%%%%%%%%
%%%%%%%%%%%%%%%%%%%%%%%%%%%%%%%%%%%%%%%%%%%%%%%%%%%%%%%%%%%%%%%%%%%%%%%%%%%%%%%%%%%%%%%%%%%%%%%%%
%%%%%%%%%%%%%%%%%%%%%%%%%%%%%%%%%%%%%%%%%%%%%%%%%%%%%%%%%%%%%%%%%%%%%%%%%%%%%%%%%%%%%%%%%%%%%%%%%
%%%%%%%%%%%%%%%%%%%%%%%%%%%%%%%%%%%%%%%%%%%%%%%%%%%%%%%%%%%%%%%%%%%%%%%%%%%%%%%%%%%%%%%%%%%%%%%%%
%%%%%%%%%%%%%%%%%%%%%%%%%%%%%%%%%%%%%%%%%%%%%%%%%%%%%%%%%%%%%%%%%%%%%%%%%%%%%%%%%%%%%%%%%%%%%%%%%
 %%%%%%%%%%%%%%%%%%%%%%%%%%%%%%%%%%%%%%%%%%%%%%%%%%%%%%%%%%%%%%%%%%%%%%%%%%%%%%%%%%%%%%%%%%%%%%%%

\vskip 0,1cm \noi  {\gsec 2010 AMS Mathematical Subject Classification}: Primary: 11M26 ; 
Secondary: 26D05, 60G17.  \par\noi  
{{\gsec Keywords and phrases}: Riemann zeta function, zeros, Dirichlet polynomials, Tur\'an's
localization criterion, increments, local suprema.

  \section{\bf 
%Introduction-
Main Result}

 \medskip  \noi       
    The question of the  existence of an  eventual explicit relation between the zeros of the Riemann Zeta function $\zeta(s)$,
 $s=\s+it$ and the prime numbers    was    raised  already  by  Landau in \cite{[L]}. 
%    The   existence of an  eventual explicit relation between the zeros of the Riemann Zeta function $
%\zeta(s)$,
% $s=\s+it$ and the prime numbers is a question which was   originally raised    by Landau  in \cite{[L]}. 
Motivated by Landau's remark, Tur\'an had much investigated the connection between zerofree regions of  $\zeta $  and local bounds of
Dirichlet polynomials expanded over the primes, see \cite{[Tu]} and 
\cite{[Tu1]}, Chapters 33-36. 
% \smallskip\par
Among the several strong localization results stated in
\cite{[Tu]}, the following semi-global criterion (Theorem 3') is of particular relevance in the present work.   
\medskip\par 
\noi{\bf Tur\'an's Localization Criterion.}  {\it Let $D$ be some positive real and $0<E\le 9/10$. Suppose there exist  positive reals   $T, \b$,
$0<\b<1$  such that for $T-T^E\le \tau\le    T+T^E$, the inequality
 \begin{equation}\label{Tu} \Big|\sum_{N_1\le p\le N_2}p^{-i\tau} \Big|\le c\  {N\log^{10} N\over \tau ^\b},
\end{equation}
holds for  
$$T^{D(1-\b^{1/6})}\le N\le N_1<N_2\le 2N\le
T^{D(1+\b^{1/6})}$$
where $c$ stands for positive numerical, explicitely calculable constant. 
 \smallskip\par   
 Then $\zeta(s)\not=0$ in the parallelogram $\s>1-\b^2 $, $T-T^E \le t\le    T+T^E$.}
 \medskip\par In this article, we show by using a local randomization argument, that Tur\'an's   approach for   localizing zeros of
$\zeta$  is sufficiently powerful to permit  to identify a completely new semi-global zerofree region. 
 
\smallskip\par
Our main result states:
%We prove the following theorem  which is our main result.
\begin{theorem} \label{tI}  Let $0<\a^*<1$. There exist $1/2<\s_0<1$, $B\ge 4$,    $\nu_0<\infty$,
such that:
\vskip 2pt
  For all
$\nu\ge
\nu_0$, there exists at least  $\a^*\, 2^{ B\nu+1}$ indices $j$ for which
 $$\zeta(\s+it)\not=0   
 \qq \forall \s\ge
\s_0, \ \ \forall t\in [2^{2B\nu} + (j-1)2^{ B\nu-1},2^{2B\nu} + j2^{ B\nu-1}[.
$$  \end{theorem}
It follows from the proof  that any value $\s_0 >1-1/(19)^{12} $ is for instance suitable.  
% it is improvable and that
     The same   approach permits to get  only     slightly better thresholds.
%\begin{theorem} \label{tI}   There exist $1/2<\s_0<1$  such that for any $0<\a<1$, there exists $
%\nu_0<\infty$, such that   $\zeta(\s +it)\not=0$ for all $\s\ge \s_o$ and all $t\in K_i$
% where $$K_i= \Big[\psi(2^{2B\nu})+ (i-1)\frac{ 2^{ B\nu}}{M}  , \psi(2^{2B\nu})+ i\frac{2^{ B\nu}}{M}   %\Big[ , \quad     1\le i\le M\big( 2^{ B\nu} +3(\sqrt 2 -1)\big), $$ 
%and $\psi(t) = t+ 3\sqrt t$,  for     at least 
% $$     \a M\big( 2^{ B\nu} +3(\sqrt 2 -1)\big) . $$
%indices $i$,  and this is true for all $\nu\ge \nu_0$, all $M\ge 2$ integer.\end{theorem}
\vskip 3pt

%  \medskip\par\noi 
   In order to bound 
$|\sum_{N_1\le p\le N_2}p^{-i\tau
 }|$, uniformly over  a family of suitable segments $[N_1, N_2]$ of the real line, we use  an approach
%  which combines an argument from   the
%stochastic processes theory with mean  value properties Dirichlet polynomials and  a variant of  Hilbert's inequality    due to  
%Montgomery and Vaughan.    
 %   This
which  can be described as follows. Let        
$ 
\p_1,
\ldots,
\p_N $ be    distinct  reals.
      Consider a finite family of  Dirichlet polynomials
   $P^s(t)=  \sum_{n=1}^N c^s_ne^{it\p_n}$, $s\in S$,   $c^s_1,\ldots, c^s_N$ being   complex numbers.     
Instead of directly searching a bound of $\sup_S|P^s(t)|$, uniformly in $t$ over some finite interval $L$, we       operate  with
 the {\it shifted} Dirichlet polynomials
  \begin{equation}P^s_\t( t)= \sum_{n=1}^N c^s_ne^{i(\t+t)\p_n},\label{not}
  \end{equation} where $\t$ belongs to   some fixed interval $J$, and $\t$ will
be treated as a random parameter. 
%   Our underlying probability  space will be $J$  equipped with the normalized Lebesgue measure.
 Given some interval $L$,      $\{P^s_\t( t),   s\in S, t\in L, \t\in J\} $  is considered at some intermediate stage of the proof, as
a random process built on $J$, of which we estimate the increments by means of variant form  of  Hilbert's inequality    due to  
 Montgomery and Vaughan.  
A classical argument from random processes machinery, allows to   efficiently    control    suprema, namely here 
 $ \sup_{t\in L} \sup_S|P^s_\t( t)|$. 
\vskip 2pt
Another step is devoted to carefully adjusting some inherent family of parameters,  in order to apply Tur\'an's result.
Once this is achieved,    a family of   intervals $(I_\t)_\t$ free of zeros is  then exhibited.   The family  is   indexed  by  a
measurable set of
$\t$'s  of controlable positive measure.  Finally, a covering argument allows  to establish  the existence of a semi-global region. 
 This is the strategy we    apply. 
 %%%%%%%%%%%%%%%%%%%%%%%%%%%%%%%%%%%%%%%%%%%%%%%%%%%%%%%%%%%%%%%%%%%%%%%%%%%%%%%%%%%%%%%%%%%%%%%%%%%%%%%%%%%%%%%%%%%
%%%%%%%%%%%%%%%%%%%%%%%%%%%%%%%%%%%%%%     %%%%%%%%%%%     %%%%%%%%%%%     %%%%%%%%%%%%%%%%%%%%%%%%%%%%%%%%%%%%%%%%
%%%%%%%%%%%%%%%%%%%%%%%%%%%%%%%%%%%%%%     %%%%%%%%%%%     %%%%%%%%%%%     %%%%%%%%%%%%%%%%%%%%%%%%%%%%%%%%%%%%%%%%
%%%%%%%%%%%%%%%%%%%%%%%%%%%%%%%%%%%%%%     %%%%%%%%%%%     %%%%%%%%%%%     %%%%%%%%%%%%%%%%%%%%%%%%%%%%%%%%%%%%%%%%
%%%%%%%%%%%%%%%%%%%%%%%%%%%%%%%%%%%%%%%%%%%%%%%%%%%%%%%%%%%%%%%%%%%%%%%%%%%%%%%%%%%%%%%%%%%%%%%%%%%%%%%%%%%%%%%%%%%
%%%%%%%%%%%%%%%%%%%%%%%%%%%%%%%%%%%%%%%%%%%%%%%%%%%%%%%%%%%%%%%%%%%%%%%%%%%%%%%%%%%%%%%%%%%%%%%%%%%%%%%%%%%%%%%%%%%
%%%%%%%%%%%%%%%%%%%%%%%%%%%%%%%%%%%%%%%%%%%%%%%%%%%%%%%%%%%%%%%%%%%%%%%%%%%%%%%%%%%%%%%%%%%%%%%%%%%%%%%%%%%%%%%%%%%
%%%%%%%%%%%%%%%%%%%%%%%%%%%%%%%%%%%%%%%%%%%%%%%%%%%%%%%%%%%%%%%%%%%%%%%%%%%%%%%%%%%%%%%%%%%%%%%%%%%%%%%%%%%%%%%%%%%  
\section{\bf Local Mean Value Results}
      Let $q$ be some positive integer and denote 
$$E_q=\Big\{\underline{   k}=(k_1,\ldots, k_N); 
k_i\in \N : k_1+\ldots+ k_N= q\Big\}.$$ Let        
$ 
\p_1,
\ldots,
\p_N $ be     linearly independent  reals. Introduce a {\it coefficient  of linear spacing} of order
$q$  by putting
 $$ \xi_\p(N, q) =\inf_{\underline{   h}, \underline{   k}\in E_q\atop 
\underline{   h}\not= \underline{   k} } \big|(h_1 -k_1)\p_1+\ldots+   (h_N-k_N)\p_N \big|.  $$
    By assumption $\xi_\p(N, q) >0$ and
 $ \xi_\p(N, 1) =\inf\{ | \p_i -\p_j  |:   
i\not=j \}  $.      In the   case  $\p_n=\log p_n$,
 $p_n$ denoting the $n$-th consecutive prime, we have the classical estimate
  $ \xi_{\p }(N, q) \ge  
  p_N^{-q}$, see before (\ref{minxi}) for a proof. 
 \vskip 3pt 
     We estimate the local increments   of $P_.$ defined in (\ref{not}). Let $J$ be a bounded
interval  and let $|J|$ denote  its length. Let $m_{J}$ denote
    the normalized Lebesgue measure on   $J$.   With the notation (\ref{not}),   if $J=[a,b]$ then
$\big\|P_.(t)-P_.(s) \big\|_{m_{J} , 2q}$ and $\big\|P_.(t)  \big\|_{m_{J} , 2q}$  respectively denote 
$$  \Big({1\over b-a}\int_{a }^{b }\big|P(\t+t)-P(\t+s)\big|^{2q}d\t \Big)^{1/2q}, \qq \Big({1\over b-a}\int_{a }^{b }\big|P(\t+t)
\big|^{2q}d\t 
\Big)^{1/2q}.$$  
  
 Introduce    the stationary metric     on the real line defined by  
$$ d(s,t)= d_N(s,t):=\Big( 2\sum_{n=1}^N
|c_n|^2\big| \sin{(   t-s) \p_n\over 2}\big|^2
\Big)^{1/ 2 }.  $$
     
\begin{proposition}    \label{p1}
   a) For any  reals   $s  $ and $t$,
$$\big\|P_.(t)-P_.(s) \big\|_{m_{J} , 2q}\le \Big( q! +{2 \min(N^{ q},\pi q!)\over |J|\xi}
\Big)^{1/2q}d(s,t)
 .  $$
 And   
$$ \big\| P_.(t) \big\|_{m_{J} , 2q}\le       \Big( q! +{2 \min(N^{
q},\pi q!)\over T\xi} \Big)^{1/2q}   \Big(\sum_{n=1}^N|c_n|^{ 2 }\Big)^{1/2} 
  . $$
 \end{proposition}
 By taking $J=[-T,T]$, $t=0$ in the last estimate, we deduce  
\begin{corollary} \label{c1}We have the following bound
$$ {1\over 2T}\int_{-T }^{T }\Big|\sum_{n=1}^N c_ne^{i \t  \p_n}
\Big|^{2q}d\t  
 \le     q!  \Big( 1 +{2 \pi   \over T\xi_\p(N, q)} \Big)   \Big(\sum_{n=1}^N|c_n|^{ 2 }\Big)^{q} 
    . $$
In particular,  
$$ {1\over 2T}\int_{-T }^{T }\Big|\sum_{p\le N }  {c_p\over    p ^{i\t}}
\Big|^{2q}d\t 
 \le     q!  \Big( 1 +{2 \pi N^q  \over T} \Big)   \Big(\sum_{p\le N}|c_p|^{ 2 }\Big)^{q} 
    . $$
\end{corollary} 
 \medskip\par
 Now  put 
$$\B=\B_\p(J,N,q)=\Big[q!  \Big( 1 +{2  \pi  \over |J|\xi_\p(N,q) } \Big)\Big]^{1/2q}. $$

\begin{theorem}\label{t2}      Let  $\tilde{
\p}_N=
\sup_{n\le N} |\p_n|$. There exists a constant
$C_q$ depending on $q$ only, such that for any interval $L$, 
 \begin{eqnarray*} \big\|\sup_{t \in L}| P_.(t)| \big\|_{m_{J },2q} & \le& C_q\, \B\max\Big\{1, |L|\tilde{ \p}_N \Big\}^{1/ 2q}
\bigg\{\Big[\sum_{n=1}^N |c_n|^2 
\Big]^{1/ 2} +  \cr &{}& \qq\qq\qq   \min\Big(|L|,  {1\over \tilde{ \p}_N} \Big)   \Big[\sum_{n=1}^N |c_n|^2 \p_n^2
\Big]^{1/ 2}\bigg\} . 
\end{eqnarray*}
\end{theorem}

%%%%%%%%%%%%%%%%%%%%%%%%%%%%%%%%%%%%%%%%%%%%%%%%%%%%%%%%%%%%%%%%%%%%%%%%%%%%%%%%%%%%%%%%%%%%%%%%%%%%%%%%%%%%%%%%%%%%%%%%%%%
%%%%%%%%%%%%%%%%%%%%%%%%%%%%%%%%%%%%%%%%%%%%%%%%%%%%%%%%%%%%%%%%%%%%%%%%%%%%%%%%%%%%%%%%%%%%%%%%%%%%%%%%%%%%%%%%%%%%%%%%%%%
\begin{proof}[Proof of Proposition \ref{p1}]   Let
$J =[d,d+T]$. Write more shortly $\xi= \xi_\p(N, q)$. 
    Plainly 
\begin{eqnarray*}\big(P(\t+t)-P(\t+s)\big)^{ q}&=&\big(\sum_{n=1}^N c_ne^{i \t  \p_n} \big( e^{i t \p_n} -e^{i s \p_n}\big) \big)^{
q}\cr&=&
\sum_{\underline{   k}\in E_q}  {q!\over k_1!\ldots k_N!} \prod_{n=1}^Nc_n^{ k_n} e^{i \t k_n\p_n } \big( e^{i t \p_n} -e^{i s
\p_n}\big)^{ k_n}
\end{eqnarray*} 
Put $\gamma_n=  e^{i t \p_n} -e^{i s
\p_n}  $. Thus
\begin{eqnarray}\label{debut} && \big|P(\t+t)-P(\t+s)\big|^{2q}\!\!  \cr&=& \!\! \sum_{\underline{   k},\underline{   h}\in E_q} 
{(q!)^2\over k_1!h_1!\ldots k_N!h_N!} 
 \prod_{n=1}^Nc_n^{ k_n}{\overline{c }_n} ^{ h_n} e^{i \t (k_n-h_n)\p_n } \gamma_n^{ k_n}  \overline{  \gamma}_n^{\, h_n} \cr & =&
\sum_{\underline{   k} \in E_q } \Big({ q!  \over k_1! \ldots k_N! } 
\Big)^2\prod_{n=1}^N(|c_n| | \gamma_n |)^{ 2k_n} +R (\t)
\end{eqnarray} 
where 
\begin{equation}\label{3}R(\t)  =  \sum_{\underline{   k},\underline{   h}\in E_q\atop 
 \underline{   k}\not=\underline{   h}} \Big({(q!)^2\over k_1!h_1!\ldots k_N!h_N!} 
\Big)\prod_{n=1}^N(c_n\gamma_n)^{ k_n}(\overline{c_n\gamma_n})^{ h_n} e^{i \t (k_n-h_n)\p_n  } .
\end{equation}  
    Owing to linear independence $\sum_{n=1}^N(k_n-h_n)\p_n=0$, iff $k_n=h_n$,   $n=1,\ldots, N$.    By integrating 
\begin{eqnarray}\label{ow}& &{1\over T}\int_J\big|P(\t+t)-P(\t+s)\big|^{2q}d\t  = \sum_{\underline{   k} \in E_q } \Big({ q! \over k_1! \ldots
k_N! } 
\Big)^2\prod_{n=1}^N    (|c_n| | \gamma_n |)^{ 2k_n} \cr
& &  + \sum_{\underline{   k},\underline{   h}\in E_q\atop 
 \underline{   k}\not=\underline{   h}}  {(q!)^2\over k_1!h_1!\ldots k_N!h_N!} 
\prod_{n=1}^N(c_n\gamma_n)^{ k_n}(\overline{c_n\gamma_n})^{ h_n}
 \cr & &\qq\qq \times\Big[{  e^{i(d+ T)
\sum_{n=1}^N(k_n-h_n)\p_n }-e^{i  d
\sum_{n=1}^N(k_n-h_n)\p_n }  \over iT(\sum_{n=1}^N(k_n-h_n)\p_n)}\Big] .
\end{eqnarray}  
 Put    
$${\bf c}_{\underline{   k}}= \prod_{n=1}^N{(c_n\gamma_n e^{i (d+T)   \p_n })^{ k_n}\over k_n!  }, \quad {\bf d}_{\underline{  
k}}= \prod_{n=1}^N{(c_n\gamma_n e^{i  d    \p_n })^{ k_n}\over k_n!  } ,\quad  {\bf   l}_{\underline{   k}}=\sum_{n=1}^N k_n \p_n .$$  
Then \begin{eqnarray}\label{4}&& {1\over T}\int_J\big|P(\t+t)-P(\t+s)\big|^{2q}d\t \cr&=& q!^2
\sum_{\underline{   k} \in E_q } |{\bf d}_{\underline{   k}}|^2  
 +{(q!)^2\over iT}\bigg\{\sum_{\underline{   k},\underline{   h}\in
E_q\atop 
 \underline{   k}\not=\underline{   h}}    { {\bf c}_{\underline{   k}}\overline{{\bf c}}_{\underline{   h}}\over  
{\bf   l}_{\underline{   k}}-{\bf   l}_{\underline{   h}}}    -\sum_{\underline{   k},\underline{   h}\in E_q\atop 
\underline{   k}\not=\underline{   h}}    { {\bf d}_{\underline{   k}}\overline{{\bf d}}_{\underline{   h}}\over   
{\bf   l}_{\underline{   k}}-{\bf   l}_{\underline{   h}}}    \bigg\} .
\end{eqnarray} 
   
 \medskip Each of the two claimed bounds will now be deduced from either Hilbert's inequality or Cauchy-Schwarz inequality.
Recall  Hilbert's inequality (\cite{[Mon]}, p.138): 
\smallskip{\it Let $\l_1,\ldots, \l_N$ be distinct real numbers, and suppose
 that $\d>0$ is chosen so that $|\l_m-\l_n|\ge \d$ whenever $n\not=m$. Then}
\begin{equation}\label{hi}\Big|\sum_{1\le m,n\le N\atop n\not=m}{x_my_n\over \l_m-\l_n}\Big|\le {\pi\over \d}\Big(\sum_{m=1}^N
|x_m|^2\Big)^{1/2}\Big(\sum_{n=1}^N |y_n|^2\Big)^{1/2}. 
\end{equation}
We shall apply it under the following form:     let $\{x_{\underline{   k} } ,   y_{\underline{   k} }, \underline{  
k}\in E_q\}
$. Let also
$\{\l_{\underline{   k} },\underline{   k} \in E_q\} $  be distinct  real numbers such that $\min\{|\l_{ \underline{   k} }-\l_{
\underline{   h} }|,{ \underline{   k}}\not={
\underline{   h} }\}\ge \d$.
 Let   $\nu =\#\{E_q\}$ and consider  a bijection   
  $i: 
 \{ 1,\ldots,\nu\}\to E_q$. 
 By  using (\ref{hi}) 
\begin{eqnarray}\label{him}\Big|\sum_{\underline{   k},\underline{   h}\in E_q\atop  \underline{   k}\not=\underline{  
h}}{x_{\underline{   k} }y_{\underline{   h} }\over
 \l_{ \underline{   k} }-\l_{
 \underline{   h} }}\Big| & =&\Big|\sum_{1\le u,v\le \nu\atop   u\not=v}{x_{i(u)
}y_{i(v) }\over
 \l_{ i(u) }-\l_{ i(v) }}\Big| \cr &\le& {\pi\over
\d}\Big(\sum_{1\le u \le \nu}  |x_{ i(u )}|^2\Big)^{1/2}\Big(\sum_{1\le  v\le \nu}  |y_{i(v) }|^2\Big)^{1/2}
\cr &= & {\pi\over
\d}\Big(\sum_{\underline{   k} \in E_q}  |x_{\underline{   k}  }|^2\Big)^{1/2}\Big(\sum_{  \underline{   h}\in E_q} 
|y_{\underline{  h} }|^2\Big)^{1/2}. 
\end{eqnarray}

By applying Hilbert's inequality to each of the two sums   in   parenthesis of the right-term  in (\ref{4}), we obtain  
\begin{equation}\label{hima} {(q!)^2\over  T}\bigg|\sum_{\underline{   k},\underline{   h}\in E_q\atop 
 \underline{   k}\not=\underline{   h}}    { {\bf c}_{\underline{   k}}\overline{{\bf c}}_{\underline{   h}}\over  
{\bf   l}_{\underline{   k}}-{\bf   l}_{\underline{   h}}}   -\sum_{\underline{   k},\underline{  
h}\in E_q\atop 
 \underline{   k}\not=\underline{   h}}    { {\bf d}_{\underline{   k}}\overline{{\bf d}}_{\underline{   h}}\over   
{\bf   l}_{\underline{   k}}-{\bf   l}_{\underline{   h}}}    \bigg|\le 
  {2\pi(q!)^2\over
  T\xi}\sum_{\underline{   k} \in E_q }      |{\bf d}_{\underline{   k}}|^2 \le  {2\pi q!\over
  T\xi}\ d(s,t)^{2q} ,
\end{equation} 
since
\begin{eqnarray}\label{8}(q!)^2\sum_{\underline{   k} \in E_q }      |{\bf d}_{\underline{   k}}|^2&=&\sum_{ k_1+\ldots+k_N=q } 
  \Big[ {  q!  \over k_1! \ldots k_N!
} \Big]^2   \prod_{n=1}^N|c_n\gamma_n|^{ 2k_n}  \cr &\le &q!\,\sum_{ k_1+\ldots+k_N=q } 
    {  q!  \over k_1! \ldots k_N!
}    \prod_{n=1}^N|c_n\gamma_n|^{ 2k_n}   = q!\,\Big[\sum_{n=1}^N |c_n\gamma_n|^2\Big]^q\cr &=& q!\Big[\sum_{n=1}^N |c_n|^2|e^{i t \p_n}
-e^{i s
\p_n}|^2\Big]^q \cr &=& q!\,\Big[4\sum_{n=1}^N |c_n|^2|\sin   {(t -s)\p_n \over 2}|^2\Big]^q
=q!\,d(s,t)^{2q} .
\end{eqnarray}
 Similarly as before
 \begin{equation}\sum_{\underline{   k} \in E_q } \Big({ q! \over k_1! \ldots
k_N! } 
\Big)^2\prod_{n=1}^N|c_n|^{ 2k_n}   \big| e^{i t \p_n} -e^{i s
\p_n}\big|^{ 2k_n}\le q!\Big[\sum_{n=1}^N |c_n\gamma_n|^2\Big]^q=q!d(s,t)^{2q}.
\end{equation}

By substituting in (\ref{4}), we therefore  get
\begin{equation}{1\over T}\int_J\big|P(\t+t)-P(\t+s)\big|^{2q}d\t \le     q!\Big(1+   {2\pi \over
  T\xi}
\Big) d(s,t)^{2q} .
\end{equation}
   Without Hilbert's inequality, it is possible to arrive to a similar result. We have with (\ref{3}), (\ref{8})  
 \begin{eqnarray*} & &{1\over T}\int_J\big|P(\t+t)-P(\t+s)\big|^{2q}d\t \le q!d(s,t)^{2q}  \cr & &  + \sum_{\underline{  
k},\underline{   h}\in E_q\atop 
 \underline{   k}\not=\underline{   h}}  {(q!)^2\over k_1!h_1!\ldots k_N!h_N!} 
\prod_{n=1}^N(c_n\gamma_n)^{ k_n}(\overline{c_n\gamma_n})^{ h_n} \cdot \Big|{  e^{i T
\sum_{n=1}^N(k_n-h_n)\p_n }-1  \over iT(\sum_{n=1}^N(k_n-h_n)\p_n)}\Big|  
\cr & \le & q!d(s,t)^{2q}+{2\over   T\xi}\Big(2 \sum_{n=1}^N  |c_n \sin{(   t-s) \p_n\over 2}|\Big)^{ q}\Big(2
\sum_{n=1}^N  |c_n
\sin{(   t-s) \p_n\over
 2}|\Big)^{ q}
 \cr & = & q!d(s,t)^{2q}+{2\over T\xi}\Big(2 \sum_{n=1}^N  |c_n \sin{(   t-s) \p_n\over 2}|\Big)^{2 q}
   \cr & \le & \Big( q! +{2 N^{ q}\over T\xi} \Big)d(s,t)^{2q}
   ,
\end{eqnarray*}
where we used Cauchy-Schwarz inequality for getting the last estimate. Combining the two last estimates gives
\begin{equation}{1\over T}\int_J\big|P(\t+t)-P(\t+s)\big|^{2q}d\t \le \Big( q! +{2 \min(N^{ q},\pi q!)\over T\xi}
\Big)d(s,t)^{2q}.
\end{equation}
 Hence the first in
assertion a). 
  The same proof also yields, mutatis mutandis
\begin{equation}  {1\over T}\int_J\big| P(\t+s )\big|^{2q}d\t \le    \big(\sum_{n=1}^N|c_n|^{ 2 }\big)^{q}  \Big( q! +{2 \min(N^{ q},\pi
q!)\over T\xi}
\Big)    .
\end{equation}
 We start with  
\begin{eqnarray*} P(\t+t) ^{ q}=\big(\sum_{n=1}^N c_ne^{i \t  \p_n}  e^{i t \p_n}  \big)^{
q}=
\sum_{\underline{   k}\in E_q}  {q!\over k_1!\ldots k_N!} \prod_{n=1}^Nc_n^{ k_n} e^{i \t k_n\p_n }  e^{i t \p_n k_n}   
\end{eqnarray*} 
and put this time $\gamma_n=  e^{i t \p_n}   $. Then all calculations made after (\ref{debut}) remain valid.
\end{proof}
 
\begin{proof}[Proof of Corollary \ref{c1}] The first assertion is immediate. As for the second, we have to estimate 
$$ \xi_\p(N, q) =\inf_{\underline{   h}, \underline{   k}\in E_q\atop 
\underline{   h}\not= \underline{   k} } \big|(h_1 -k_1)\p_1+\ldots+   (h_N-k_N)\p_N \big|.  $$
when  $\p_n=\log p_n$. Let $\underline{   \ell}=\underline{   h} -\underline{   k}$ and 
put 
$$    P^+=  \prod_{ \ell_n> 0} p_n^{\ell_n}, \qq P^-=  \prod_{\ell_n< 0} p_n^{-\ell_n}$$
Notice that $ P^+\not= P^-$ by assumption, and $\max (P^+, P^-) \le p_N^q$. Suppose $P^+> P^- $. Then
$$\big|\ell_1\p_1+\ldots+  \ell_N\p_N \big|= \big|\log \prod_{  n =1}^Np_n^{\ell_n}\big|   = \log{ P^+\over P^-}
% =  \log1 +{ P^+-P^-\over P^-}
\ge   \log1 +{
1\over P^-}\ge   \log1 +{
1\over p_N^q}\ge {
1\over p_N^q}.
$$   
  The case  $P^+< P^-$ is   treated identically. Therefore  \begin{equation}\label{minxi}
   \xi_\p(N,q)   \ge      p_N^{-q} 
   .
   \end{equation}
 And so, it suffices to apply the first estimate to this case.
\end{proof}%%%%%%%%%%%%%%%%%%%%%%%%%%%%%%%%%%%%%%%%%%%%%%%%%%%%%%%%%%%%%%%%%%%%%%%%%%%%%%%%%%%%%%%%%%%%%%%%%%%%%%%%%%%%%%%%%%%%%%%%%%%%%%%%%%%%%%%
%%%%%%%%%%%%%%%%%%%%%%%%%%%%%%%%%%%%%%%%%%%%%%%%%%%%%%%%%%%%%%%%%%%%%%%%%%%%%%%%%%%%%%%%%%%%%%%%%%%%%%%%%%%%%%%%%%%%%%%%%%%%%%%%%%%%%%%
%%%%%%%%%%%%%%%%%%%%%%%%%%%%%%%%%%%%%%%%%%%%%%%%%%%%%%%%%%%%%%%%%%%%%%%%%%%%%%%%%%%%%%%%%%%%%%%%%%%%%%%%%%%%%%%%%%%%%%%%%%%%%%%%%%%%%%%
%%%%%%%%%%%%%%%%%%%%%%%%%%%%%%%%%%%%%%%%%%%%%%%%%%%%%%%%%%%%%%%%%%%%%%%%%%%%%%%%%%%%%%%%%%%%%%%%%%%%%%%%%%%%%%%%%%%%%%%%%%%%%%%%%%%%%%%
%%%%%%%%%%%%%%%%%%%%%%%%%%%%%%%%%%%%%%%%%%%%%%%%%%%%%%%%%%%%%%%%%%%%%%%%%%%%%%%%%%%%%%%%%%%%%%%%%%%%%%%%%%%%%%%%%%%%%%%%%%%%%%%%%%%%%%%
%%%%%%%%%%%%%%%%%%%%%%%%%%%%%%%%%%%%%%%%%%%%%%%%%%%%%%%%%%%%%%%%%%%%%%%%%%%%%%%%%%%%%%%%%%%%%%%%%%%%%%%%%%%%%%%%%%%%%%%%%%%%%%%%%%%%%%%

\begin{proof}[Proof of Theorem \ref{t2}] We  
need some elements from the theory of stochastic processes. See \cite{[W1]},  also \cite{[W4]} and references therein for a similar treatment.  Let
$(T,\d)$ be a compact  metric space and denote by $D$ the diameter of $T$. For any $x\in T$ and $\e>0$, let $B(x,\e)$ denote the  open $\d$-ball of $T$ with
center $x$ and radius
$\e$.
  A stochastic process
$X=\{X(t),t\in T\}$ is simply a collection of random variables indexed by $T$, and defined on some common probability space $(\O, {\mathcal
A},\P)$. Let
$1<p<\infty$. Consider the increment condition
\begin{equation} \|X(s)-X(t)\|_p\le d(s,t)\qquad (s,t\in T)\label{2.15}\end{equation}
  
 Assume that there exists a probability measure
$\m$ on
$T$ such that 
\begin{equation} \sup_{x\in T}\int_0^D {d\e \over \m(B(x,\e))^{1/p}}= M<\infty .\label{2.16}
\end{equation} 
 By Theorem 4.6 in \cite{[Ta]},  each separable process that satisfies the increment condition (\ref{2.15}), is sample
continuous. Moreover 
\begin{equation}\big\|\sup_{s,t\in T} |X(s)-X(t)| \big\|_p \le K_p M, \label{2.17}\end{equation}
where $K_p$ depends on $p$ only. The above inequality   follows  from
the majorizing measure condition   (\ref{2.16}) and Proposition 2.7  in \cite{[Ta]}. The sample continuity property is in turn obtained
by   combining Theorem 4.6 with Theorem 2.9   in  \cite{[Ta]}. A stochastic process is separable (with respect to  $\d$), if there
exists a countable  dense subset
$T_0$ of $T$ such that for each $t$ in $T$, 
 $ X(t) =\lim_{T_0\ni s\rightarrow t  }X(s)$, almost surely.  
 By Proposition \ref{p1}    
 $$\big\|P_.(t)-P_.(s) \big\|_{m_{J} , 2q}\le  \B\, d(s,t), \qq  \big\| P_.(s) \big\|_{m_{J} , 2q}\le       \B  \,
\Big(\sum_{n=1}^N|c_n|^{ 2 }\Big)^{1/2} 
  . $$
 The trajectories $s\mapsto
P_\t (s)$ being continuous for   every  $\t$,   $P_.$ is thus trivially separable.  
   As
     $d ^2(s,t)  \le 
 4\pi^2 |s-t |^2\sum_{n=1}^N |c_n|^2\big(   \p_n^2\wedge {1\over \pi^2 |s-t|^2}\big)$,  
 we   deduce that 
  \begin{equation}  d (s,t)
\le 2\pi  |s-t| \big(\sum_{n=1}^N |c_n|^2 \p_n^2 \big)^{1/2},\label{2.18}\end{equation} 
once $\pi |s-t| \le 1/\tilde{ \p}_N$.
   Consider a covering $\{I_{j },j=j_1,\ldots, j_1+H\}$    of $L$ with intervals 
 $$ I_{ j} =[{j-1\over \pi\tilde{ \p}_N } , {j\over \pi\tilde{ \p}_N}[ \qq\qq   (j\ge 1)  .$$  
  Introduce an auxiliary process ${\mathcal Y}$ defined for $s \in I_{ j}$, $j\ge 1$ by
$$  {\mathcal Y}_s  = { P_.(s)-P_.({j-1\over \pi\tilde{ p}_N})  \over 2\pi    \B \big(\sum_{n=1}^N |c_n|^2 \p_n^2 \big)^{1/2}}         
    .$$ 
 By  (\ref{2.18}), for every $s,t\in I_{ j}$
\begin{equation} \| {\mathcal Y}_s-{\mathcal Y}_t  \|_{m_{J },2q}={\|  P_.(s)-P_.(t)  \|_{m_{J },2q} \over 2\pi    \B\big(\sum_{n=1}^N |c_n|^2
\p_n^2
\big)^{1/2}}  \le {d(s,t) \over 2\pi \B \big(\sum_{n=1}^N |c_n|^2 \p_n^2 \big)^{1/2}}  \le
 |s-t| .\label{2.19}\end{equation} 

Thus $\{ {\mathcal Y}_s,s \in I_{ j}\}$ satisfies  (\ref{2.15})  with the usual metric. 
Recall that $m_{I_j}$ denotes the normalized Lebesgue measure on $I_j$. Then
$$ \int_0^{ {\rm diam}(I_{ j} )  }{d\e\over  m_{I_j} (B(s,\e)) ^{1/  2q  }}\le \int_0^{1/  (\pi\tilde{ p}_N)  } \Big({1\over 
\pi\tilde{ p}_N\e  }\Big)^{ 1/  2q   }d\e =
%\buildrel{\e={\eta\over  \pi\tilde{ p}_N}}\over {=}
 {1\over \pi\tilde{ p}_N} \int_0^{1  } 
\eta ^{  {-1/  2q  }}d\e
 \le     {c_q\over  \tilde{ p}_N} .$$

 Hence$$\sup_{s\in I_{ j}}\int_0^{{\rm diam}(I_{ j} )}{d\e\over  m_{I_j} (B(s,\e)) ^{1/ 2q }} \le    {c_q\over  \tilde{ p}_N}  .$$
 From (\ref{2.17}) follows that   
\begin{equation} \sup_{j=1}^N\big\|\sup_{s,t\in I_{ j}}
\left|{\mathcal Y}_s-{\mathcal Y}_t\right|\big\|_{m_{J },2q}  \le   {  c'_q\over  \tilde{ p}_N}. 
\label{2.20}\end{equation}  

 Assume that $|L|\pi \tilde{ \p}_N> 1 $, and let $\{I_{j },j=j_1,\ldots, j_1+H\}$, $H\ge 0$, be a covering of $L$.
 Let $s\in L$, and let $j$ be such that $s\in I_j$. By writing 
$$P_.(s)=P_.( {j-1\over \pi\tilde{ p}_N})+\big(P_.(s)-P_.(  {j-1\over \pi\tilde{ p}_N} )\big)=P_.( {j-1\over
\pi\tilde{ p}_N})+2\pi   \B  \big(\sum_{n=1}^N |c_n|^2 \p_n^2 \big)^{1/2} {\mathcal Y}_s ,$$
 next using the triangle
inequality, we get  
\begin{eqnarray}\big\|\sup_{s \in L}| P_.(s)| \big\|_{m_{J },2q} &\le &\big\|\sup_{1\le j\le H}|P_.( {j-1\over \pi\tilde{ p}_N})|
\big\|_{m_{J },2q}  \cr && + 2\pi    \B \big(\sum_{n=1}^N |c_n|^2 \p_n^2 \big)^{1/2}\big\|\sup_{1\le j\le H\atop s \in I_{ j}}|{\mathcal
Y}_s|
\big\|_{m_{J },2q}.\label{2.21}\end{eqnarray}

In the one hand  
\begin{equation}\big\|\sup_{1\le j\le H}|P_.( {j-1\over \pi\tilde{ p}_N})| \big\|_{m_{J },2q}\le  H^{1\over 2q}\sup_{1\le j\le H}\big\| 
P_.( {j-1\over \pi\tilde{ p}_N})  \big\|_{m_{J },2q}\le  \B H^{1\over 2q}\Big(\sum_{n=1}^N |c_n|^2  \Big)^{1\over2}
.\label{2.22}\end{equation}  And in the other
\begin{equation}\big\|\sup_{1\le j\le H\atop s \in I_{ j}}|{\mathcal Y}_s| \big\|_{m_{J },2q}\le H^{1\over 2q}\sup_{1\le j\le H 
}\big\|\sup_{  s \in I_{ j}}|{\mathcal Y}_s|
\big\|_{m_{J },2q}  .\label{2.23}\end{equation}
 But ${\mathcal Y}({j-1\over \pi\tilde{ \p}_N})=0$, and so by  (\ref{2.20})
 $$ \|\sup_{s\in I_{ j}}
 |{\mathcal Y}_s | \|_{m_{J },2q} \le  \|\sup_{s,t\in I_{ j}}
\left|{\mathcal Y}_s-{\mathcal Y}_t\right| \|_{m_{J },2q}  \le  {  c'_q\over  \tilde{ p}_N} .$$    
   
 As $H \le C \max(1, |L|\tilde{ \p}_N )$, we deduce 
% (letting $C _q= \max(1,c'_q)$)
\begin{equation} \big\|\sup_{s \in L}| P_.(s)| \big\|_{m_{J },2q}  \le C_q\B \,  (|L|\tilde{ \p}_N  )^{1\over
2q}\bigg\{\Big(\sum_{n=1}^N |c_n|^2 
\Big)^{1/2} +     {1\over \tilde{ \p}_N }   \big(\sum_{n=1}^N |c_n|^2 \p_n^2 \big)^{1/2}\bigg\} .\label{2.24}\end{equation}

 Finally, if $|L|\pi \tilde{ \p}_N\le 1 $, write $L=[L_1, L_2]$.  Given $s,t\in L$, we have $\pi |s-t| \le \pi |L| \le 1/\tilde{ \p}_N$,
and so 
 $$ \|P_.(s)-P_.(t)\||_{m_{J },2q}\le  \B \,d (s,t)
\le 2\pi   \B |s-t| \big(\sum_{n=1}^N |c_n|^2 \p_n^2 \big)^{1/2}.$$  
Put 
$${\mathcal P}_s={  P_.(s)- P_.(L_1)\over 2\pi   \B  \big(\sum_{n=1}^N |c_n|^2 \p_n^2 \big)^{1/2}}, \qq s\in L. $$  
Then
 $\| {\mathcal P}_s-{\mathcal P}_t\|_{m_{J },2q}\le |s-t|$. 
Similarly as for getting  (\ref{2.20}), we obtain
\begin{equation}  \big\|\sup_{s,t\in L}
\left|{\mathcal P}_s-{\mathcal P}_t\right|\big\|_{m_{J },2q}  \le  c_q|L|. 
\label{2.25}\end{equation}
It follows that 
\begin{equation} \big\|\sup_{s \in L}| P_.(s)| \big\|_{m_{J },2q}  \le C'_q \B\,\bigg\{\big(\sum_{n=1}^N |c_n|^2 
\big)^{1/2} +     |L|  \big(\sum_{n=1}^N |c_n|^2 \p_n^2 \big)^{1/2}\bigg\} .\label{2.37}\end{equation}
With  (\ref{2.24}) and  (\ref{2.16}), we arrived to
 \begin{eqnarray} \big\|\sup_{s \in L}| P_.(s)| \big\|_{m_{J },2q} & \le &C''_q\B\,\max\big\{1, |L|\tilde{ \p}_N \big\}^{1\over 2q}
\Big\{\big[\sum_{n=1}^N |c_n|^2 
\big]^{1\over 2} \cr & &+     \min\Big(|L|,  {1\over \tilde{ \p}_N} \Big)   \Big[\sum_{n=1}^N |c_n|^2 \p_n^2 \Big]^{1\over
2}\Big\} . \label{2.26}\end{eqnarray}
  This achieves the proof.
   \end{proof}

%%%%%%%%%%%%%%%%%%%%%%%%%%%%%%%%%%%%%%%%%%%%%%%%%%%%%%%%%%%%%%%%%%%%%%%%%%%%%%%%%%%%%%%%%%%%%%%%%%%%%%%%%%%%%%%%%%%%%%%%%%%%%%%%%%%%%%%
%%%%%%%%%%%%%%%%%%%%%%%%%%%%%%%%%%%%%%%%%%%%%%%%%%%%%%%%%%%%%%%%%%%%%%%%%%%%%%%%%%%%%%%%%%%%%%%%%%%%%%%%%%%%%%%%%%%%%%%%%%%%%%%%%%%%%%%
%%%%%%%%%%%%%%%%%%%%%%%%%%%%%%%%%%%%%%%%%%%%%%%%%%%%%%%%%%%%%%%%%%%%%%%%%%%%%%%%%%%%%%%%%%%%%%%%%%%%%%%%%%%%%%%%%%%%%%%%%%%%%%%%%%%%%%%
%%%%%%%%%%%%%%%%%%%%%%%%%%%%%%%%%%%%%%%%%%%%%%%%%%%%%%%%%%%%%%%%%%%%%%%%%%%%%%%%%%%%%%%%%%%%%%%%%%%%%%%%%%%%%%%%%%%%%%%%%%%%%%%%%%%%%%%
%%%%%%%%%%%%%%%%%%%%%%%%%%%%%%%%%%%%%%%%%%%%%%%%%%%%%%%%%%%%%%%%%%%%%%%%%%%%%%%%%%%%%%%%%%%%%%%%%%%%%%%%%%%%%%%%%%%%%%%%%%%%%%%%%%%%%%%
%%%%%%%%%%%%%%%%%%%%%%%%%%%%%%%%%%%%%%%%%%%%%%%%%%%%%%%%%%%%%%%%%%%%%%%%%%%%%%%%%%%%%%%%%%%%%%%%%%%%%%%%%%%%%%%%%%%%%%%%%%%%%%%%%%%%%%%
\section{\bf Proof of Theorem \ref{tI}}
  The constants appearing in Tur\'an's result  (Section 1) are important. We have therefore explicited  all constants appearing in our proof. 
  %It will be  also of matter to   clarify when a constant is numerical   because of the numerical constant %$c$ in (\ref{Tu}), whose effective value has no interest. 
 \vskip 3pt 
 We begin with applying Theorem \ref{t2}  to 
$$P(N_1, N_2,t)=\sum_{N_1\le p\le N_2}p^{-it  } $$ where $N\le   N_1<N_2\le 2N $.  We have $ 
\tilde{
\p}_N\le
\sup\{  \log p, p\le 2N\}  \le C\log N$ and by using (\ref{minxi}),
\begin{equation}\B\le \Big(q!\big[  1 +{2  \pi p_{N_1}^q \over |J|  } \big]\Big)^{1/2 q}
 \le C_q\max\big(1,{   N^q  \over |J|  }
\big)^{1/2 q} .\label{3.1}\end{equation}
  Let 
$L $ be such that $|L| \ge 1$. Since  $\pi(2x)-\pi(x)\le {x\over \log x}$ for any integer $x>1$, we have $\pi(N_2)-\pi(N_1)\le \pi(2N 
)-\pi(N  )<  N /\log N
$, we have 
$$ \sum_{N  \le p\le 2N  } \log^2 p\le  \log^2 (2N)\sum_{N  \le p\le 2N  }1\le  \frac{N \log^2 (2N)}{\log N} \le CN \log N.$$
  We    get
\begin{eqnarray}\Big\|\sup_{t \in L}\big |P_.(N_1, N_2,t)\big |\Big\|_{m_{J },2q }     & \le&
C_q\max\big(1,{   N^q  \over |J|  }
\big)^{1/ 2 q}   (|L|\log N   )^{1/ 2q }   
\bigg\{\Big({N \over \log N }\Big)^{1/ 2} 
\cr & &\qq     +    
   {1\over \log  N  }\Big(   \sum_{N  \le p\le 2N  } \log^2 p \Big) ^{1/ 2 }\bigg\} . 
\cr & \le&
C_q  \Big(\max\big(1,{   N^q  \over |J|  }
\big)       |L|\log N   \Big)^{1/ 2q }  \Big({N \over \log N }\Big)^{1/ 2}     . 
\end{eqnarray}
So that if $|J|\le N^q$,  
\begin{equation} \big\|\sup_{t \in L}\big |P_.(N_1, N_2,t)\big |\big\|_{m_{J },2q }     \le
C_q {N \over (\log N)^{1/ 2} }     \Big({|L|\log N   \over |J|} \Big)^{1/ 2q }       . \label{3.2}\end{equation}
The remainding part of the proof  now consists  of carefully adjusting the parameters in order  to apply
Tur\'an's result (\ref{Tu}).  

 \medskip\par
\noi {\bf Main parameters: ($\bf H,\d, q, B, \nu, m, \a$).} The constants $H,\d, q, \a$ are   numerical  and fixed. They will produce the constant $c$ in
(\ref{Tu}). See (\ref{3.12a}). 
  \vskip 3pt 
  
  Let $H\ge 2$ be some  integer. Put
 $$ \d= {H-1\over 8 H}\qq q={5\over  1 -8\d  }=5H .$$
  Then 
$$0<\d<1/8 \qq {\rm and}  \qq q> {4(\d+1)\over  1 -8\d  }.$$
%The precise admissible range of values of $H$ is given by 
% \begin{equation}\label{H} 3\le H\le \frac{8^7}{ 5}e^{-5/3} .
%  \end{equation}   
In addition  we set   
 $$B=  4q \d   +2(\d+1) ,$$
  and notice that 
  $2B=  8q \d   +4(\d+1)<q $. 
  \vskip 3pt 
  Now  fix some positive integer  $\nu$  and set
   $$U=2^\nu,  \qq J= [U^{2B}, 2U^{2B}] , \qq  L= [U^{ B}, 8U^{ B}] .  $$ 
 Let $N=2^m$  with $m\ge \nu$.  It follows that   
 $|J|=U^{2B}\le U^{q} \le N^q$.    Then
\begin{equation}\Big\|\sup_{2^m\le   N_1<N_2\le 2^{m+2}}\sup_{t \in L}\big |P_.(N_1, N_2,t)\big |\Big\|_{m_{J},2q} \le C_q{2^{ m(1+1/ q)
}    \over m^{1/2}} \Big({|L|  m   \over |J|}
\Big)^{1/ 2q }.\label{3.3}\end{equation}
   
By Minkowski's inequality
\begin{eqnarray*}&&\Big\|\sup_{\nu \le m\le \nu(1+\d) }\sup_{  2^m\le  
N_1<N_2\le 2^{m+2}}\sup_{t \in L}\big |P_.(N_1, N_2,t)\big |\Big\|_{m_{J},2q}\cr
&\le &\Big\|\sum_{\nu \le m\le \nu(1+\d) }\sup_{  2^m\le  
N_1<N_2\le 2^{m+2}}\ \sup_{t \in L}\big |P_.(N_1, N_2,t)\big |\Big\|_{m_{J},2q}
\cr& \le& C_q \Big({|L|   
\over |J|}
\Big)^{1/ 2q } \sum_{\nu \le m\le
\nu(1+\d)}2^{ m(1+1/ q) }    
     m  
   ^{1/ 2q -1/2} 
\cr& \le& C_q \nu  
 ^{1/ 2q -1/2}\Big({|L|   
\over |J|}
\Big)^{1/ 2q } \sum_{\nu \le m\le
\nu(1+\d)}2^{ m(1+1/ q) }    
 \cr& \le& 2C_q \nu   
 ^{1/ 2q-1/2 } 2^{-(B/ 2q)\nu}  2^{\nu(1+\d)(1+1/ q) }   .    
\end{eqnarray*}
   Now  if   $U\le N\le N_1<N_2\le 2N\le U^{1+\d}$, choose $\nu \le m\le
\nu(1+\d)$   such  that $2^m\le N< 2^{m+1} $. Then $2^m\le N\le N_1<N_2\le 2N< 2^{m+2} $. Thus
\begin{eqnarray}&&\Big\|\sup_{U\le N\le N_1<N_2\le 2N\le U^{1+\d}  }\ \sup_{  t \in
L}\big |P_.(N_1, N_2,t)\big |\Big\|_{m_{J},2q}\cr &\le&
\Big\|\sup_{\nu
\le m\le \nu(1+\d) }\ \sup_{  2^m\le   N_1<N_2\le 2^{m+2}}\ \sup_{t \in L}
\big |P_.(N_1, N_2,t)\big |\Big\|_{m_{J},2q}\cr &\le& 2C_q    2^{\nu[(1+\d)(1+1/ q) -(B/ 2q)]}  \nu   
 ^{1/ 2q-1/2 }    \cr &\le& 2C_q    2^{ [1-\d]\nu}  \nu   
 ^{1/ 2q-1/2 }   := M .
\end{eqnarray}
  since with our choices  $(1+\d)(1+1/ q)-B/2q=1-\d$.    
  
  \vskip 5 pt
   Next let  $0<\a<1$ be fixed   and set  $\m(\a)=1/(1-\a)^{ 1/(2q)}$. Set 
   $$\tilde J=  \Big\{\t\in J: \sup_{U\le N\le N_1<N_2\le 2N\le U^{1+\d} \atop t \in L}\big |P_\t(N_1, N_2,t)\big |\le \m(\a)
  M\Big\} .$$
   By the Tchebycheff inequality
\begin{eqnarray} {1\over |J|}\l \{J\backslash \tilde J  \} &\le& {1\over |J| (\m M)^{  2q}}\int_J
\sup_{U\le N\le N_1<N_2\le 2N\le U^{1+\d} \atop t \in L}\big |P_\t(N_1, N_2,t)\big |^{  2q} d\t 
\cr & \le  &\m(\a)  ^{- 2q}= 1-\a .
\end{eqnarray}
Therefore  $\l \{\tilde J  \}\ge \a|J|$ and for all $\t\in \tilde J$,
\begin{equation}\sup_{U\le N\le N_1<N_2\le 2N\le U^{1+\d} \atop t \in L}\big |P_\t(N_1, N_2,t)\big |\le  
     2 \m(\a)C_q  2^{ [1-\d]\nu}  \nu   
 ^{1/ 2q-1/2 }  . \label{3.6}
\end{equation}

Pick some $\t$ in $\tilde J$. Then  
\begin{eqnarray}\sup_{U\le N\le N_1<N_2\le 2N\le U^{1+\d} \atop \tau  \in \t +L}\big |\sum_{N_1\le p\le N_2}{1\over p^{i\tau}}\big
|&\le  &  2 \m(\a)C_q  2^{\nu(1-\d) }   \nu  ^{1/ 2q-1/2 }\cr&=&2 \m(\a)C_q    U^{ 
1-\d }   (\log U)  ^{1/ 2q-1/2 }
\cr &\le& 2 \m(\a)C_q    \ {U   (\log U)  
^{1/ 2q -1/2}\over U^{   \d  }}. 
\end{eqnarray}  
 But if $\tau  \in \t +L$, $\tau\le 2 U^{2B}+8U^B \le 3U^{2B}$  if $U$, namely $\nu$ is large enough. It follows that $U^\d\ge C
\tau^{
\d/ (2B)}$. 

\vskip 3 pt Put  
$$b:={\d \over 2B}={\d \over 8q \d   +4(\d+1)} . $$
We have obtained:
\vskip 3pt 
For all $\tau \in [\t + U^B, \t+8U^B] $ and $U\le N\le N_1<N_2\le 2N\le U^{1+\d}$,   
\begin{equation} \big |\sum_{N_1\le p\le N_2}{1\over p^{i\tau}}\big |\le  2 \m(\a)C_q   \ {N   (\log N)  
^{1/ 2q-1/2 }\over \tau^b}.\label{3.8}
\end{equation} 
 \medskip\par
\noi \noi {\bf A family of local zerofree regions:}  
We use secondary  parameters: $\d_0, D, b$.  Let 
$$T=T_\t= \t +3\sqrt \t . $$We may assume $ \t\ge 1$. In the one hand
$$T-\sqrt T=\t +3\sqrt \t-\sqrt \t\sqrt{1 +{3/ \sqrt \t}}\ge \t +3\sqrt \t- 2\sqrt \t =\t + \sqrt \t\ge \t+ U^B .$$
And in the other since $U^{2B}\le \t\le 2U^{2B}$
$$T+\sqrt T=\t +3\sqrt \t+\sqrt \t\sqrt{1 +{3/ \sqrt \t}}\le \t +5\sqrt \t \le \t +5\sqrt2 U^B \le \t +8  U^B .$$
Hence $[T-\sqrt T, T+\sqrt T]\subset \t +L$ and   estimate (\ref{3.8}) is valid for $T-\sqrt T\le   \tau \le T+\sqrt T$. Further, as  
$$U^{2B}\le \t\le T=\t +3\sqrt \t\le 2U^{2B}+3\sqrt 2 U^{ B}=U^{2B}[2+ 3\sqrt 2 U^{- B}]\le 7  U^{2 B},$$  
it is also valid in the restricted range of values 
\begin{equation}\label{3.9}T^{1\over 2B}\le N\le N_1<N_2\le 2N\le \Big({T\over 7}\Big)^{1+\d\over 2B}. \end{equation}

 Now select a positive real $\d_0$ such
that  
$$0<{2\d_0\over 1-\d_0}<\d 
 .$$
  We notice that $1+\d- {1+\d_0\over 1-\d_0}=\d- {2\d_0\over 1-\d_0}>0$. Choose $ \nu$  sufficiently large so that $2^{  \nu [\d-
{2\d_0\over 1-\d_0}]}\ge 7^{1+\d}$. Since $2B>1$  we have 
% $$T ^{ 1+\d- {1+\d_0\over 1-\d_0} }=T ^{ \d- {2\d_0\over 1-\d_0} }\ge 2^{2B(\d- {2\d_0\over 1-\d_0})\nu}\ge 2^{ \nu [1+\d- {1+\d_0\over
%1-\d_0}]}\ge 7^{1+\d},$$
$$T ^{ 1+\d- {1+\d_0\over 1-\d_0} } \ge 2^{2B\nu ( 1+\d- {1+\d_0\over 1-\d_0}) }\ge 2^{ \nu [1+\d- {1+\d_0\over
1-\d_0}]}=2^{ \nu [\d- {2\d_0\over 1-\d_0}]}\ge 7^{1+\d},$$ namely
$$\Big({T\over 7}\Big)^{1+\d}  \ge  T ^{   {1+\d_0\over 1-\d_0} }.$$  

Next  put  
$$ D={1\over 2B(1-\d_0)}.$$ 
Then   (\ref{3.9})  implies the admissibility of the more suitable field of parameters
\begin{equation}T^{D(1-\d_0)}=T^{1\over 2B}\le N\le N_1<N_2\le 2N\le T ^{ D(1+\d_0)}=T^{ {1+\d_0\over 2B( 1-\d_0)}}\le \Big({T\over
7}\Big)^{1+\d\over 2B}.\label{3.10}\end{equation} 
 Estimate (\ref{3.8}) then implies 
 \begin{equation} \big |\sum_{N_1\le p\le N_2}{1\over p^{i\tau}}\big |\le  2 \m(\a)C_q   \ {N   (\log N)  
^{1/ 2q-1/2 }\over \tau^b},\label{3.11}\end{equation}
for all $\tau \in [  T-T^{1/2},  T+T^{1/2}] $ and all $T^{D(1-\d_0)}\le N\le N_1<N_2\le 2N\le T^{D(1+\d_0)}$.
 \vskip 3pt 
Recall that $0<\d<1/8$  and $q={5\over  1 -8\d  }$. Thus 
$$B=  4q \d   +2(\d+1)<  {20\d\over   1 -8\d   }    +{9\over 4}= {80 \d  +9-72\d\over  4(1 -8\d)  } =  {8 \d +9 \over  4(1 -8\d)  
}<  {5\over  2(1 -8\d ) }.$$ And $$b={\d \over 2B}\ge {\d( 1 -8\d ) \over 5} . $$
In order that $b^{1/6}\ge   \d_0 $, it suffices that ${\d( 1 -8\d ) \over 5}\ge (\d/2)^6 $, namely $   1 -8\d    \ge (5/2^6) \d  ^5$, 
%or
%since $  1 -8\d \ge 1/9$, that $  \d  ^5\le {2^6 \over  90}       $. But
%$$  \d  ^5\le {1\over 9^5}<{2^6 \over  90} .      $$
which is fulfilled if $\d<1/9$ for instance, namely recalling that $\d= \frac{H-1}{8H}$, if  $H<9$ which we do. 
%$$ 3\le H\le \frac{8^7}{ 5}e^{-5/3} .$$
%  which is  precisely ensured by  (\ref{H}). 
  
  Thus $b\ge   \d_0^6$ does hold, and (\ref{3.11}) implies  that the inequality
 \begin{equation} \big |\sum_{N_1\le p\le N_2}{1\over p^{i\tau}}\big |\le  c   \ {N   (\log N)  
 ^{1/ 2q-1/2 }\over \tau^{\d_0^6}},
 \label{3.12}\end{equation}
 %\begin{equation}\sup_{T-T^{1/2}\le \tau \le   T+T^{1/2}\atop T^{D(1-\d_0)}\le N\le
%N_1<N_2\le 2N\le T^{D(1+\d_0)}} \big |\sum_{N_1\le p\le N_2}{1\over p^{i\tau}}\big |\le  c   \ %{N   (\log N)  
%^{1/ 2q-1/2 }\over \tau^{\d_0^6}},
%\label{3.12}\end{equation}
with  (recalling that   $\m(\a)=1/(1-\a)^{ 1/(2q)}$),
 \begin{equation}\label{3.12a} c= 2 \m(\a)C_q,
\end{equation}
  \vskip 3pt \noi  holds for all $\tau \in [  T-T^{1/2},  T+T^{1/2}] $ and all $T^{D(1-\d_0)}\le N\le
N_1<N_2\le 2N\le T^{D(1+\d_0)}$.

  \medskip\par
  Tur\'an's result (section 1) then implies   that 
    \begin{equation} \label{3.13}\zeta(\s+it)\not= 0, \qq \forall \s>1-\d_0^{12}, \  \  \forall t\in [T_\t-T_\t^{1/2},   T_\t+T_\t^{1/2}].\end{equation}

But this holds for {\it any} $\t\in \tilde J$ (recalling that $\l (\tilde J)\ge \a|J|   $, $J= [2^{2B\nu},  2^{2B\nu +1}]$), and for    {\it any}   $\nu$, assuming this one   large enough,
depending on
$\d$, say $\nu_\d$. We also recall that $\d$ was fixed from the beginning (see "Main parameters").
\begin{remark} Finding  {\it one} $\t$ in $J$ such that $\zeta(\s+it)\not=0$ for all $t$ in  $[T_\t-T_\t^{1/2},   T_\t+T_\t^{1/2}]$ and $\s>\s_0$, for some $\s_0<1$, can be deduced from   Carlson's estimate on the number of zeros of the Riemann zeta function.
%Notice that for a given fixed $\t$, (\ref{3.13}) follows from Carlson's result on the number of zeros of the %Riemann zeta function.
 The point here is that  we have    a measurable set of values of $\t$'s   of measure close to the one of $J$,  for which this is valid. This together with 
a simple covering argument will permit to exhibit a much bigger zerofree zone.
\end{remark}

  \medskip\par
\noi {\bf A semi-global  zerofree region:}  
%We first draw some simple facts from the validity of property (\ref{3.13}) for the family  of parameters 
%described above, on the set  $J= [2^{2B\nu},  2^{2B\nu +1}]$, 
    Let    $\psi(\t) = \t + 3  \sqrt \t$. The indice $\nu $ with   $\nu \ge \nu_\d$ being now temporarily fixed,    let  $J_0= ]2^{2B\nu},  2^{2B\nu +1}[\backslash \tilde J$. Using the fact that  
  $ \l( \psi ( [a,b])) = (b-a) +   3(\sqrt b -\sqrt a) \le (b-a) \{ 1+  2. 2^{  -B\nu   }\}$, one can show 
  \begin{equation}\label{approx}\l( \psi (J_0) ) 
      \le  \{ 1+ 1 /{ 2^{ B\nu  }}\}  (1-\a)    \l (J)    . 
\end{equation}      
    Let $\eta>0$, $J_0$ being an open set,   $J_0=\cup_{n=1}^\infty I_n $  where $I_n$ are open intervals, 
% (approximating by open sets would have suffice since $\l$ is a Radon measure). 
  Let  $U_N= \cup_{n=1}^N I_n $.  
   Writing $U=U_N{\cup  \!\!\!\cdot\ } B$   with
$B\subset  \cup_{n=N+1}^\infty I_n$,
we have,
\begin{eqnarray*} \l( \psi (J_0) )&\le &  \l\big( \psi (U_N){\cup  \!\!\!\cdot\ } \psi (B) \big)   \le \l( \psi (U_N)) +  \sum_{n=N+1}^\infty \l(   \psi
(   I_n))
\cr  &\le & \l( \psi (U_N)) +  \{ 1+2. 2^{  -B\nu   }\} \sum_{n=N+1}^\infty \l(       I_n )
     \le     \l( \psi (U_N)) +  \eta \{ 1+2. 2^{  -B\nu   }\}  ,\end{eqnarray*}
 assuming $N$ large enough.   
   Further $\cup_{n=1}^N I_n={\cup  \!\!\!\cdot\ }_{n=1}^{N'} I'_n$,     $I'_n$ being pairwise disjoint intervals.  Since $\psi $ is continuous increasing,
  $$ \l( \psi (U_N))  =\l( \sum_{n=1}^{N'}  \psi (I'_n))=\sum_{n=1}^{N'} \l(    \psi (I'_n)) \le  \{ 1+ 2. 2^{  -B\nu   }\} \sum_{n=1}^{N } \l(     I_n )
  $$ $$= \{ 1+ 2. 2^{  -B\nu   }\}\l(  U_N )\le   \{ 1+ 2. 2^{  -B\nu   }\}(\l(J_0)+  \eta). $$
   
Thus  
$$  \l( \psi (J_0) ) 
      \le   \{ 1+ 1 /{ 2^{ B\nu  }}\} \l(J_0) +2\eta \{ 1+ 1 /{ 2^{ B\nu  }}\}    
     \le  \{ 1+ 1 /{ 2^{ B\nu  }}\} \big\{ (1-\a)    \l (J) +2\eta\big\} ,$$
since $\l (  J_0)\le (1-\a)\l (J)   $.  Since $\eta$ is arbitrary,    (\ref{approx}) follows. 
 \vskip 3pt  
 Therefore,  \begin{eqnarray} \label{3.14}  \l( \psi (  \tilde J) ) &\ge& \l( \psi (J ) )- \{ 1+    2^{ -B\nu   }\}  (1-\a)    \l (J)
\cr &=& \l( \psi (J ) )\Big[1 - \frac{  1+  { 2^{ -B\nu  }}  }{1+ 3(\sqrt 2 -1) 2^{- B\nu}} (1-\a)   \Big]
  \cr &:=&  (1-\bar \a) \l( \psi (J ) ),
         \end{eqnarray} 
noticing that $ \l( \psi (    J) ) 
%=2^{2B\nu}+ 3(\sqrt 2 -1)2^{ B\nu} 
=\l(     J  )(1+ 3(\sqrt 2 -1) 2^{- B\nu} ) $. 
 \vskip 5pt 
    As $T_\t^{1/2}\ge \t^{1/2}\ge 2^{ B\nu}$, we have   $[T_\t-T_\t^{1/2},   T_\t+T_\t^{1/2}]\supset [T_\t- 2^{ B\nu},   T_\t+2^{ B\nu}]$.
 Now consider on  $\psi (    J)=  [\psi(2^{2B\nu}),  \psi(2^{2B\nu +1})]  $   the subdivision 
$$K_i= \Big[\psi(2^{2B\nu})+ (i-1)   2^{ B\nu -1}   , \psi(2^{2B\nu})+ i2^{ B\nu -1}   \Big[ , \quad     1\le i\le  \big( 2^{ B\nu +1} +6(\sqrt 2
-1)\big). $$ 
   In view of (\ref{3.14}),  the number of indices $i$ such that $K_i\cap \psi(\tilde J))= \emptyset$ is less than 
$ (1-\bar \a) \l(\psi (    J)) /2^{ B\nu +1}$. 
\vskip 5pt   Consequently, at least $  \bar \a  \l(\psi (    J)) /2^{ B\nu+1}$ indices $i$ are such that $K_i\cap \psi(\tilde J))\not = \emptyset$. Pick a real   
$\vartheta $ in  the intersection. We have  
$$[\vartheta-\vartheta^{1/2},   \vartheta+\vartheta^{1/2}]\supset K_i.$$ 
  
 \vskip 3pt So that by (\ref{3.13}),
   \begin{equation} \label{3.15}\zeta(\s+it)\not= 0, \qq \forall \s>1-\d_0^{12}, \  \  \forall t \in  K_i, \end{equation}
and the number of indices $i$ for  which  this is true,   exceeds   
 \begin{equation} \label{3.16}  \bar \a  \l(\psi (    J)) /2^{ B\nu+1}=\bar \a  \big( 2^{ B\nu+1} +6(\sqrt 2 -1)\big) . 
 \end{equation}

  \vskip 3 pt
We can now achieve the proof. Given any fixed real $0<\a^*<1$,  it follows from (\ref{3.15}),(\ref{3.16}) that  in any subdivision of $\psi (    J)$ of size $  2^{ B\nu-1} $,
  at least $   \a^*  2^{ B\nu+1} $
  intervals are free of zero. 
Since $\psi (    J)=[2^{2B\nu}+3.2^{ B\nu},2.2^{2B\nu }+3\sqrt 2.2^{ B\nu}]$, it also implies that  in any subdivision of $[2^{2B\nu},2^{2B\nu +1[} $ of
size
$  2^{ B\nu-1} $,
  at least $   \a^*    2^{ B\nu+1}  $ intervals are free of zero.

\bigskip
 
 {
}

 \end{document}